# Asymptotic expansions at any time for scalar fractional SDEs with Hurst index $H > 1/2$

SÉBASTIEN DARSES[1] and IVAN NOURDIN[2]

[1]*Boston University, Department of Mathematics and Statistics, 111 Cummington Street, Boston, MA 02215, USA. E-mail: darses@math.bu.edu*

[2]*Université Pierre et Marie Curie, Laboratoire de Probabilités et Modèles Aléatoires, 4 place Jussieu, Boîte courrier 188, 75252 Paris Cedex 5, France. E-mail: ivan.nourdin@upmc.fr*

We study the asymptotic expansions with respect to $h$ of

$$\mathrm{E}[\Delta_h f(X_t)], \qquad \mathrm{E}[\Delta_h f(X_t)|\mathscr{F}_t^X] \quad \text{and} \quad \mathrm{E}[\Delta_h f(X_t)|X_t],$$

where $\Delta_h f(X_t) = f(X_{t+h}) - f(X_t)$, when $f:\mathbb{R} \to \mathbb{R}$ is a smooth real function, $t \geq 0$ is a fixed time, $X$ is the solution of a one-dimensional stochastic differential equation driven by a fractional Brownian motion with Hurst index $H > 1/2$ and $\mathscr{F}^X$ is its natural filtration.

*Keywords:* asymptotic expansion; fractional Brownian motion; Malliavin calculus; stochastic differential equation

## 1. Introduction

We study the asymptotic expansions with respect to $h$ of

$$\begin{aligned}
P_t f(h) &\triangleq \mathrm{E}[\Delta_h f(X_t)], \\
\widehat{P}_t f(h) &\triangleq \mathrm{E}[\Delta_h f(X_t)|\mathscr{F}_t^X], \\
\widetilde{P}_t f(h) &\triangleq \mathrm{E}[\Delta_h f(X_t)|X_t],
\end{aligned} \qquad (1)$$

with $\Delta_h f(X_t) \triangleq f(X_{t+h}) - f(X_t)$, when $f:\mathbb{R} \to \mathbb{R}$ is a smooth real function, $t \geq 0$ is a fixed time, $X$ is the solution to the fractional stochastic differential equation

$$X_t = x + \int_0^t b(X_s)\,\mathrm{d}s + \int_0^t \sigma(X_s)\,\mathrm{d}B_s, \qquad t \in [0,T], \qquad (2)$$

and $\mathscr{F}^X$ is its natural filtration. Here, $b, \sigma : \mathbb{R} \to \mathbb{R}$ are real functions belonging to the space $\mathrm{C}_b^\infty$ of all bounded continuous functions having bounded derivatives of all order, while $B$ is a one-dimensional fractional Brownian motion with Hurst index $H \in (1/2, 1)$. When the integral with respect to $B$ is understood in the Young sense, equation (2) has a unique pathwise solution $X$ in the set of processes whose paths are Hölder continuous of index $\alpha \in (1 - H, H)$. Moreover, by, for example, [12], Theorem 4.3, we have, for any $g : \mathbb{R} \to \mathbb{R} \in C_b^\infty$,

$$g(X_t) = g(x) + \int_0^t g'(X_s)\sigma(X_s)\,\mathrm{d}B_s + \int_0^t g'(X_s)b(X_s)\,\mathrm{d}s, \qquad t \in [0, T]. \qquad (3)$$

The asymptotic expansion of $\mathrm{E}[f(X_h)]$ with respect to $h$ has been recently studied in [1], [7]. In our framework, it turns out to be the case where $t = 0$ since we obviously have

$$\mathrm{E}[f(X_h)] - f(x) = P_0 f(h) = \widehat{P}_0 f(h) = \widetilde{P}_0 f(h).$$

In these last references, the authors work in a multidimensional setting and under the weaker assumption that the Hurst index $H$ of the fractional Brownian motion $B$ is greater than $1/3$ (the integral with respect to $B$ is then understood in the rough paths sense of Lyons' type for [1] and of Gubinelli's type for [7]). In particular, it is proved in [1], [7] that there exists a family

$$\Gamma = \{\Gamma_{2kH+\ell} : (k, \ell) \in \mathbb{N}^2, (k, \ell) \neq (0, 0)\}$$

of differential operators such that for any smooth $f : \mathbb{R} \to \mathbb{R}$, we have the following asymptotic expansion:

$$P_0 f(h) \underset{h \to 0}{\sim} \sum h^{2kH+\ell} \Gamma_{2kH+\ell}(f, \sigma, b)(x). \qquad (4)$$

Moreover, in [7], operators $\Gamma_{2kH+\ell}$ are expressed using trees.

A natural question now arises. Can we also get an expansion of $P_t f(h)$ when $t \neq 0$? Let us first consider the case where $B$ is the standard Brownian motion (which corresponds to the case where $H = 1/2$). By the Markov property on one hand, we have $\widehat{P}_t f(h) = \widetilde{P}_t f(h)$. On the other hand, we always have $P_t f(h) = \mathrm{E}[\widehat{P}_t f(h)]$. Thus, there exist relations between $P_t f(h)$, $\widehat{P}_t f(h)$ and $\widetilde{P}_t f(h)$. Moreover, the asymptotic expansion of $P_t f(h)$ can be obtained as a corollary to that of $P_0 f(h)$ using the conditional expectation either with respect to the past $\mathscr{F}_t^X$ of $X$, or with respect to $X_t$ only and the strong Markov property.

When $H > 1/2$, $B$ is not Markovian. The situation regarding $P_t f(h)$, $\widehat{P}_t f(h)$ and $\widetilde{P}_t f(h)$ is then completely different and actually more complicated. In particular, we no longer have $\widehat{P}_t f(h) = \widetilde{P}_t f(h)$ and we cannot deduce the asymptotic expansion of $P_t f(h)$ from that of $P_0 f(h)$.

The current paper is concerned with the study of possible asymptotic expansions of the various quantities $P_t f(h)$, $\widehat{P}_t f(h)$ and $\widetilde{P}_t f(h)$ when $H > 1/2$. We will see that some nontrivial phenomena appear. More precisely, we will show in Section 3 that $\widehat{P}_t f(h)$ does not admit an asymptotic expansion in the scale of the fractional powers of $h$ when $t \neq 0$. Regarding $\widehat{P}_t f(h)$, the situations when $t = 0$ and $t > 0$ are thus really different. On the other hand, unlike $\widehat{P}_t f(h)$, the quantities $P_t f(h)$ and $\widetilde{P}_t f(h)$ admit, when $t \neq 0$, an asymptotic expansion in the scale of the fractional powers of $h$. However, the computation of this expansion is more difficult than in the case where $t = 0$ (as carried out in [1], [7]). That is why we prefer to consider only the one-dimensional case. As an illustration, let us consider the trivial equation $dX_t = dB_t$, $t \in [0, T]$, $X_0 = 0$. That is, $X_t = B_t$ for every $t \in [0, T]$. We have, thanks to a Taylor expansion,

$$\widetilde{P}_0 f(h) = \sum_{k=1}^{n} \frac{f^{(k)}(0)}{k!} \mathrm{E}[(B_h)^k] + \cdots = \sum_{k=1}^{\lfloor n/2 \rfloor} \frac{f^{(2k)}(0)}{2^k k!} h^{2Hk} + \cdots,$$

while by a linear Gaussian regression, when $t \neq 0$,

$$\widetilde{P}_t f(h) = \mathrm{E}\left[ f\left( \left(1 + \frac{H}{t}h - \frac{h^{2H}}{2t^{2H}} + \cdots \right) B_t + (h^{2H} + \cdots) N \right) - f(B_t) \Big| B_t \right]$$

$$= \frac{H B_t f'(B_t)}{t} h - \frac{B_t f'(B_t)}{2t^{2H}} h^{2H} + \cdots$$

with $N \sim \mathcal{N}(0, 1)$ a random variable independent of $B_t$.

One of the key points of our strategy relies on the use of a Girsanov transformation and the Malliavin calculus for fractional Brownian motion. We refer to [3], [10] for a deep insight of this topic.

We will restrict the exposition of our asymptotic expansions to the case when $\sigma = 1$. Indeed, under the assumption

(A)　　the function $\sigma$ is elliptic on $\mathbb{R}$, that is, it satisfies $\inf_{\mathbb{R}} |\sigma| > 0$

and using the change of variable formula (3), equation (2) can be reduced to a diffusion $Y$ with a constant diffusion coefficient,

$$Y_t = \int_0^{X_t} \frac{dz}{\sigma(z)}.$$

Moreover, since $\int_0^{\cdot} \frac{dz}{\sigma(z)}$ is strictly monotone from $\mathbb{R}$ to $\mathbb{R}$ under assumption (A), the $\sigma$-fields generated by $X_t$ (resp. by $X_s$, $s \leq t$) and $Y_t$ (resp. by $Y_s$, $s \leq t$) are the same. Consequently, assuming that $\sigma = 1$ is not at all restrictive since it allows the recovery of the general case under assumption (A). We therefore consider in the sequel that $X$ is the unique solution of

$$X_t = x + \int_0^t b(X_s) \, ds + B_t, \qquad t \in [0, T], \tag{5}$$

with $b \in C_b^{\infty}$ and $x \in \mathbb{R}$.

The paper is organized as follows. In Section 2, we recall some basic facts about fractional Brownian motion, the Malliavin calculus and fractional stochastic differential equations. In Section 3, we prove that $\widehat{P}_t f(h)$ does not admit an asymptotic expansion with respect to the scale of fractional powers of $h$, up to order $n \in \mathbb{N}$. We eventually show, in Section 4, that $\widetilde{P}_t f(h)$ admits an asymptotic expansion.

## 2. Preliminaries

We first briefly recall some basic facts about stochastic calculus with respect to a fractional Brownian motion. One may refer to [9], [10] for further details. Let $B = (B_t)_{t \in [0,T]}$ be a fractional Brownian motion with Hurst parameter $H \in (1/2, 1)$ defined on a probability space $(\Omega, \mathscr{A}, \mathbf{P})$. We mean that $B$ is a centered Gaussian process with the covariance function $\mathrm{E}(B_s B_t) = R_H(s,t)$, where

$$R_H(s,t) = \tfrac{1}{2}(t^{2H} + s^{2H} - |t-s|^{2H}). \tag{6}$$

We denote by $\mathcal{E}$ the set of step $\mathbb{R}$-valued functions on $[0,T]$. Let $\mathfrak{H}$ be the Hilbert space defined as the closure of $\mathcal{E}$ with respect to the scalar product

$$\langle \mathbf{1}_{[0,t]}, \mathbf{1}_{[0,s]} \rangle_{\mathfrak{H}} = R_H(t,s).$$

We denote by $|\cdot|_{\mathfrak{H}}$ the associate norm. The mapping $\mathbf{1}_{[0,t]} \mapsto B_t$ can be extended to an isometry between $\mathfrak{H}$ and the Gaussian space $\mathcal{H}_1(B)$ associated with $B$. We denote this isometry $\varphi \mapsto B(\varphi)$.

The covariance kernel $R_H(t,s)$ introduced in (6) can be written as

$$R_H(t,s) = \int_0^{s \wedge t} K_H(s,u) K_H(t,u) \, du,$$

where $K_H(t,s)$ is the square-integrable kernel defined, for $s < t$, by

$$K_H(t,s) = c_H s^{1/2 - H} \int_s^t (u-s)^{H-3/2} u^{H-1/2} \, du \tag{7}$$

with $c_H^2 = \frac{H(2H-1)}{\beta(2-2H, H-1/2)}$ and $\beta$ the Beta function. By convention, we set $K_H(t,s) = 0$ if $s \geq t$.

We define the operator $\mathcal{K}_H$ on $\mathrm{L}^2([0,T])$ by

$$(\mathcal{K}_H h)(t) = \int_0^t K_H(t,s) h(s) \, ds.$$

Let $\mathcal{K}_H^* : \mathcal{E} \to \mathrm{L}^2([0,T])$ be the linear operator defined by

$$\mathcal{K}_H^*(\mathbf{1}_{[0,t]}) = K_H(t, \cdot).$$

The following equality holds for any $\phi, \psi \in \mathcal{E}$:

$$\langle \phi, \psi \rangle_{\mathfrak{H}} = \langle \mathcal{K}_H^* \phi, \mathcal{K}_H^* \psi \rangle_{\mathrm{L}^2([0,T])} = \mathrm{E}(B(\phi) B(\psi)).$$

$\mathcal{K}_H^*$ then provides an isometry between the Hilbert space $\mathfrak{H}$ and a closed subspace of $L^2([0,T])$.

The process $W = (W_t)_{t \in [0,T]}$ defined by

$$W_t = B((\mathcal{K}_H^*)^{-1}(\mathbf{1}_{[0,t]})) \tag{8}$$

is a Wiener process and the process $B$ has the following integral representation:

$$B_t = \int_0^t K_H(t,s) \, dW_s.$$

Hence, for any $\phi \in \mathfrak{H}$, we have

$$B(\phi) = W(\mathcal{K}_H^* \phi).$$

If $b, \sigma \in C_b^\infty$, then (2) admits a unique solution $X$ in the set of processes whose paths are Hölder continuous of index $\alpha \in (1-H, H)$. Moreover, $X$ has the Doss–Sussman's-type representation (see, e.g., [5])

$$X_t = \phi(A_t, B_t), \qquad t \in [0, T], \tag{9}$$

with $\phi$ and $A$ given, respectively, by

$$\frac{\partial \phi}{\partial x_2}(x_1, x_2) = \sigma(\phi(x_1, x_2)), \qquad \phi(x_1, 0) = x_1, \qquad x_1, x_2 \in \mathbb{R}$$

and

$$A_t' = \exp\left(-\int_0^{B_t} \sigma'(\phi(A_t, s)) \, ds\right) b(\phi(A_t, B_t)), \qquad A_0 = x_0, \qquad t \in [0, T].$$

Let $b \in C_b^\infty$ and $X$ be the solution of (5). Following [11], the fractional version of the Girsanov theorem applies and ensures that $X$ is a fractional Brownian motion with Hurst parameter $H$ under the new probability $\mathbf{Q}$ defined by $d\mathbf{Q} = \eta^{-1} \, d\mathbf{P}$, where

$$\eta = \exp\left(\int_0^T \left(\mathcal{K}_H^{-1} \int_0^\cdot b(X_r) \, dr\right)(s) \, dW_s + \frac{1}{2} \int_0^T \left(\mathcal{K}_H^{-1} \int_0^\cdot b(X_r) \, dr\right)^2(s) \, ds\right). \tag{10}$$

Let $\mathscr{S}$ be the set of all smooth cylindrical random variables, that is, of the form $F = f(B(\phi_1), \ldots, B(\phi_n))$, where $n \geq 1$, $f : \mathbb{R}^n \to \mathbb{R}$ is a smooth function with compact support and $\phi_i \in \mathfrak{H}$. The Malliavin derivative of $F$ with respect to $B$ is the element of $L^2(\Omega, \mathfrak{H})$ defined by

$$D_s F = \sum_{i=1}^n \frac{\partial f}{\partial x_i}(B(\phi_1), \ldots, B(\phi_n)) \phi_i(s), \qquad s \in [0, T].$$

In particular, $D_s B_t = \mathbf{1}_{[0,t]}(s)$. As usual, $\mathbb{D}^{1,2}$ denotes the closure of the set of smooth random variables with respect to the norm

$$\|F\|_{1,2}^2 = \mathrm{E}[F^2] + \mathrm{E}[|D.F|_{\mathfrak{H}}^2].$$

The Malliavin derivative $D$ verifies the following chain rule. If $\varphi : \mathbb{R}^n \to \mathbb{R}$ is $C_b^1$ and if $(F_i)_{i=1,\ldots,n}$ is a sequence of elements of $\mathbb{D}^{1,2}$, then $\varphi(F_1,\ldots,F_n) \in \mathbb{D}^{1,2}$ and we have, for any $s \in [0,T]$,

$$D_s \varphi(F_1,\ldots,F_n) = \sum_{i=1}^{n} \frac{\partial \varphi}{\partial x_i}(F_1,\ldots,F_n) D_s F_i.$$

The divergence operator $\delta$ is the adjoint of the derivative operator $D$. If a random variable $u \in \mathrm{L}^2(\Omega, \mathfrak{H})$ belongs to the domain of the divergence operator, that is, if there exists $c_u > 0$ such that

$$|\mathrm{E}\langle DF, u\rangle_{\mathfrak{H}}| \leq c_u \|F\|_{\mathrm{L}^2} \qquad \text{for any } F \in \mathcal{S},$$

then $\delta(u)$ is defined by the duality relationship

$$\mathrm{E}(F\delta(u)) = \mathrm{E}\langle DF, u\rangle_{\mathfrak{H}}$$

for all $F \in \mathbb{D}^{1,2}$.

## 3. Study of the asymptotic expansion of $\widehat{P}_t f(h)$

Recall that $\widehat{P}_t f(h)$ is defined by (1), where $X$ is given by (5).

**Definition 1.** *We say that $\widehat{P}_t f(h)$ admits an asymptotic expansion with respect to the scale of fractional powers of $h$, up to order $n \in \mathbb{N}$, if there exist some real numbers $0 < \alpha_1 < \cdots < \alpha_n$ and some random variables $C_1,\ldots,C_n \in \mathrm{L}^2(\Omega, \mathscr{F}_t^X)$, not identically-zero, such that*

$$\widehat{P}_t f(h) = C_1 h^{\alpha_1} + \cdots + C_n h^{\alpha_n} + o(h^{\alpha_n}) \qquad \text{as } h \to 0,$$

*where $o(h^\alpha)$ stands for a random variable of the form $h^\alpha \phi_h$, with $\mathrm{E}[\phi_h^2] \to 0$ as $h \to 0$.*

If $\widehat{P}_t f(h)$ admits an asymptotic expansion in the sense of Definition 1, we must, in particular, have the existence of $\alpha > 0$ verifying the following condition:

$$\lim_{h \to 0} h^{-\alpha} \widehat{P}_t f(h) \text{ exists in } \mathrm{L}^2(\Omega) \text{ and is not identically zero.}$$

However, we have the following.

**Theorem 1.** *Let $f:\mathbb{R}\to\mathbb{R}\in C_b^\infty$ and $t\in(0,T]$. Assume, moreover, that*

$$\mathrm{Leb}(\{x\in\mathbb{R}:f'(x)=0\})=0. \tag{11}$$

*Then, as $h\to 0$, $h^{-\alpha}\widehat{P}_t f(h)$ converges in $L^2(\Omega)$ if and only if $\alpha<H$. In this case, the limit is zero.*

**Remark 1.** Since $\widehat{P}_0 f(h)=\widetilde{P}_0 f(h)$, we refer to Theorem 2 for the case where $t=0$.

**Proof of Theorem 1.** The proof is divided into two cases.

(i) *First case*: $\alpha\in(0,1]$. Since $H>1/2$, let us first remark that $h^{-\alpha}\widehat{P}_t f(h)$ converges in $L^2(\Omega)$ if and only if $h^{-\alpha}f'(X_t)\mathrm{E}[X_{t+h}-X_t|\mathscr{F}_t^X]$ converges in $L^2(\Omega)$. Indeed, we use a Taylor expansion:

$$|f(X_{t+h})-f(X_t)-f'(X_t)(X_{t+h}-X_t)|\leq \tfrac{1}{2}|f''|_\infty|X_{t+h}-X_t|^2,$$

so

$$|\widehat{P}_t f(h)-f'(X_t)\mathrm{E}[X_{t+h}-X_t|\mathscr{F}_t^X]|\leq \tfrac{1}{2}|f''|_\infty \mathrm{E}[|X_{t+h}-X_t|^2|\mathscr{F}_t^X].$$

Thus, applying Jensen's formula:

$$h^{-2\alpha}\mathrm{E}[|\widehat{P}_t f(h)-f'(X_t)\mathrm{E}[X_{t+h}-X_t|\mathscr{F}_t^X]|^2]$$
$$\leq \tfrac{1}{4}|f''|_\infty^2 h^{-2\alpha}\mathrm{E}[\mathrm{E}[|X_{t+h}-X_t|^2|\mathscr{F}_t^X]]^2$$
$$\leq \tfrac{1}{4}|f''|_\infty^2 h^{-2\alpha}\mathrm{E}[|X_{t+h}-X_t|^4]=O(h^{4H-2\alpha}).$$

Since $\alpha\leq 1<2H$, we can conclude.

By (11) and the fact that $X_t$ has a positive density on $\mathbb{R}$ (see, e.g., [8] Theorem A), we have that $h^{-\alpha}f'(X_t)\mathrm{E}[X_{t+h}-X_t|\mathscr{F}_t^X]$ converges in $L^2(\Omega)$ if and only if $h^{-\alpha}\mathrm{E}[X_{t+h}-X_t|\mathscr{F}_t^X]$ converges in $L^2(\Omega)$.

For $X$ given by (5), we have that $\mathscr{F}^X=\mathscr{F}^B$. Indeed, one inclusion is obvious, while the other can be proven using (9). Moreover, since $b\in C_b^\infty$ and $\alpha\leq 1$, the term $h^{-\alpha}\mathrm{E}[\int_t^{t+h}b(X_s)\,\mathrm{d}s|\mathscr{F}_t^X]$ converges in $L^2(\Omega)$ when $h\downarrow 0$. Therefore, we have, due to (5) that $h^{-\alpha}\mathrm{E}[X_{t+h}-X_t|\mathscr{F}_t^X]$ converges in $L^2(\Omega)$ if and only if $h^{-\alpha}\mathrm{E}[B_{t+h}-B_t|\mathscr{F}_t^B]$ converges in $L^2(\Omega)$.

Set

$$Z_h^{(t)}=h^{-\alpha}\mathrm{E}[B_{t+h}-B_t|\mathscr{F}_t^B]=h^{-\alpha}\int_0^t (K_H(t+h,s)-K_H(t,s))\,\mathrm{d}W_s,$$

where the kernel $K_H$ is given by (7) and the Wiener process $W$ is defined by (8). We have

$$\mathrm{Var}(Z_h^{(t)})=h^{-2\alpha}\int_0^t (K_H(t+h,s)-K_H(t,s))^2\,\mathrm{d}s$$

$$= h^{-2\alpha} c_H^2 \int_0^t s^{1-2H} \left( \int_t^{t+h} (u-s)^{H-3/2} u^{H-1/2} \, du \right)^2 ds.$$

We deduce

$$\begin{aligned}
\operatorname{Var}(Z_h^{(t)}) &\geq h^{-2\alpha} \left( \frac{c_H}{H-1/2} \right)^2 t^{2H-1} \\
&\quad \times \int_0^t s^{1-2H} ((t+h-s)^{H-1/2} - (t-s)^{H-1/2})^2 \, ds \\
&= h^{-2\alpha} \left( \frac{c_H}{H-1/2} \right)^2 \int_0^t \left(1 - \frac{s}{t}\right)^{1-2H} ((s+h)^{H-1/2} - s^{H-1/2})^2 \, ds \\
&= h^{2(H-\alpha)} \left( \frac{c_H}{H-1/2} \right)^2 \int_0^{t/h} \left(1 - \frac{hs}{t}\right)^{1-2H} g^2(s) \, ds
\end{aligned} \quad (12)$$

with $g(s) = (s+1)^{H-1/2} - s^{H-1/2}$. Similarly,

$$\begin{aligned}
\operatorname{Var}(Z_h^{(t)}) &\leq h^{-2\alpha} \left( \frac{c_H}{H-1/2} \right)^2 (t+h)^{2H-1} \\
&\quad \times \int_0^t s^{1-2H} ((t+h-s)^{H-1/2} - (t-s)^{H-1/2})^2 \, ds \\
&= h^{-2\alpha} \left(1 + \frac{h}{t}\right)^{2H-1} \left( \frac{c_H}{H-1/2} \right)^2 \\
&\quad \times \int_0^t \left(1 - \frac{s}{t}\right)^{1-2H} ((s+h)^{H-1/2} - s^{H-1/2})^2 \, ds \\
&= h^{2(H-\alpha)} \left(1 + \frac{h}{t}\right)^{2H-1} \left( \frac{c_H}{H-1/2} \right)^2 \int_0^{t/h} \left(1 - \frac{hs}{t}\right)^{1-2H} g^2(s) \, ds.
\end{aligned} \quad (13)$$

Note that $g^2(s) \sim (H - \frac{1}{2})^2 s^{2H-3}$ as $s \to +\infty$. So, $sg^2(s) \longrightarrow 0$ as $s \to +\infty$ and $\int_0^{+\infty} |g^2(s)| \, ds < +\infty$ since $2H - 3 < -1$. Since $s \mapsto sg^2(s)$ is bounded on $\mathbb{R}^+$, we have, by the dominated convergence theorem, that

$$\int_0^{t/h} \left( \left(1 - \frac{hs}{t}\right)^{1-2H} - 1 \right) g^2(s) \, ds = \int_0^1 \frac{(1-u)^{1-2H} - 1}{u} g^2\left(\frac{tu}{h}\right) \frac{tu}{h} \, du$$

tends to zero as $h \to 0$. Thus,

$$\lim_{h \to 0} \int_0^{t/h} \left(1 - \frac{hs}{t}\right)^{1-2H} g^2(s) \, ds = \int_0^\infty g^2(s) \, ds < +\infty.$$

Combined with (12)–(13), we now deduce that

$$\mathrm{Var}(Z_h^{(t)}) \sim h^{2(H-\alpha)} \left(\frac{c_H}{H-1/2}\right)^2 \int_0^\infty g^2(s)\,\mathrm{d}s \qquad \text{as } h \to 0. \tag{14}$$

If $Z_h^{(t)}$ converges in $\mathrm{L}^2(\Omega)$ as $h \to 0$, then $\lim_{h\to 0}\mathrm{Var}(Z_h^{(t)})$ exists and is finite. But, thanks to (14), we have that $\lim_{h\to 0}\mathrm{Var}(Z_h^{(t)}) = +\infty$ when $\alpha > H$. Consequently, $Z_h^{(t)}$ does not converge in $\mathrm{L}^2(\Omega)$ as $h \to 0$ when $\alpha > H$.

Conversely, when $\alpha < H$, we have, from (14), that $\lim_{h\to 0}\mathrm{Var}(Z_h^{(t)}) = 0$. Then $Z_h^{(t)} \xrightarrow{\mathrm{L}^2} 0$ when $\alpha < H$.

In order to complete the proof of the first case, it remains to consider the critical case where $\alpha = H$. We first deduce from (14) that $Z_h^{(t)} \xrightarrow{\mathrm{Law}} \mathcal{N}(0, \sigma_H^2)$, as $h \to 0$, with

$$\sigma_H^2 = \left(\frac{c_H}{H-1/2}\right)^2 \int_0^\infty g^2(s)\,\mathrm{d}s.$$

Let us finally show that the previous limit does not hold in $\mathrm{L}^2$. Assume for a moment that $Z_h^{(t)}$ converges in $\mathrm{L}^2(\Omega)$ as $h \to 0$. In particular, $\{Z_h^{(t)}\}_{h>0}$ is Cauchy in $\mathrm{L}^2(\Omega)$. So, denoting by $Z^{(t)}$ the limit in $\mathrm{L}^2(\Omega)$, we have $\mathrm{E}[Z_\varepsilon^{(t)} Z_\delta^{(t)}] \to \mathrm{E}[|Z^{(t)}|^2]$ when $\varepsilon, \delta \to 0$. But, for any fixed $x > 0$, we can show by using the same transformations as above that as $h \to 0$,

$$\mathrm{E}(Z_{hx}^{(t)} Z_{h/x}^{(t)}) \longrightarrow \left(\frac{c_H}{H-1/2}\right)^2 r(x) = \mathrm{E}(|Z^{(t)}|^2),$$

where

$$r(x) = \int_0^\infty ((s+x)^{H-1/2} - s^{H-1/2})\left(\left(s + \frac{1}{x}\right)^{H-1/2} - s^{H-1/2}\right)\mathrm{d}s$$
$$= x \int_0^\infty g(x^2 u) g(u)\,\mathrm{d}u.$$

Consequently, the function $r$ is constant on $]0, +\infty[$. The Cauchy–Schwarz inequality yields

$$|g|_{\mathrm{L}^2}^2 = r(1) = r(\sqrt{2}) = \langle \sqrt{2} g(2\cdot), g\rangle_{\mathrm{L}^2} \leq \sqrt{2}|g(2\cdot)|_{\mathrm{L}^2}|g|_{\mathrm{L}^2} = |g|_{\mathrm{L}^2}^2.$$

We thus have an equality in the previous inequality. We deduce that there exists $\lambda \in \mathbb{R}$ such that $g(2u) = \lambda g(u)$ for all $u \geq 0$. Since $g(0) = 1$, we have $\lambda = 1$. Consequently, for any $u \geq 0$ and any integer $n$, we get

$$g(u) = g\left(\frac{u}{2^n}\right) \xrightarrow[n\to\infty]{} g(0) = 1,$$

which is absurd. Therefore, when $\alpha = H$, $Z_h^{(t)}$ does not converge in $\mathrm{L}^2(\Omega)$ as $h \to 0$. This concludes the proof of the first case.

(ii) Second case: $\alpha \in (1, +\infty)$. If $h^{-\alpha}\widehat{P}_t f(h)$ converges in $L^2(\Omega)$, then $h^{-1}\widehat{P}_t f(h)$ converges in $L^2(\Omega)$ toward zero. This contradicts the first case, which concludes the proof of Theorem 1. □

## 4. Study of the asymptotic expansion of $\widetilde{P}_t f(h)$

Recall that $\widetilde{P}_t f(h)$ is defined by (1), where $X$ is given by (5). The main result of this section is the first point of the following theorem.

**Theorem 2.** *Let $t \in [0,T]$ and $f : \mathbb{R} \to \mathbb{R} \in C_b^\infty$. We write $\mathcal{N}$ for $\mathbb{N}^2 \setminus \{(0,0)\}$. For $(p,q) \in \mathcal{N}$, set*

$$J_{2pH+q} = \{(m,n) \in \mathcal{N} : 2mH + n \leq 2pH + q\}.$$

1. *If $t \neq 0$, there exists a family $\{Z^{(t)}_{2mH+n}\}_{(m,n)\in\mathcal{N}}$ of random variables measurable with respect to $X_t$ such that for any $(p,q) \in \mathcal{N}$,*

$$\widetilde{P}_t f(h) = \sum_{(m,n) \in J_{2pH+q}} Z^{(t)}_{2mH+n} h^{2mH+n} + o(h^{2pH+q}). \tag{15}$$

2. *If $t = 0$, for any $(p,q) \in \mathcal{N}$, we have*

$$P_0 f(h) = \widetilde{P}_0 f(h) = \widehat{P}_0 f(h)$$
$$= \sum_{(m,n)\in J_{2pH+q}} \left( \sum_{I \in \{0,1\}^{2m+n}, |I|=2m} c_I \Gamma_I(f,b)(x) \right) h^{2mH+n} + o(h^{2pH+q})$$

*with $c_I$ and $\Gamma_I$ defined, respectively, by (17) and (19) below.*

***Remark 2.*** 1. In (15), $Z^{(t)}_{0H+1}$ coincides with the stochastic derivative of $X$ with respect to its present $t$, as defined in [2].

2. The expansion (15) allows the expansion of $P_t f(h)$ to obtained for $t \neq 0$:

$$P_t f(h) = \mathrm{E}[\widetilde{P}_t f(h)] = \sum_{(m,n)\in J_{2pH+q}} \mathrm{E}[Z^{(t)}_{2mH+n}] h^{2mH+n} + o(h^{2pH+q}).$$

The following subsections are devoted to the proof of Theorem 2. Note that a quicker proof of the first assertion seems to be as follows. Once $t > 0$ is fixed, we write

$$X_{t+h} = X_t + \int_t^{t+h} b(X_s)\,ds + \widetilde{B}^{(t)}_h, \qquad h \geq 0, \tag{16}$$

where $\widetilde{B}_h^{(t)} = B_{t+h} - B_t$ is again a fractional Brownian motion. We could then think that an expansion for $\widetilde{P}_t f(h)$ directly follows from the one for $\widetilde{P}_0 f(h)$, simply by a shift. This is unfortunately not the case, due to the fact that the initial value in (16) is not just a *real number* as in the case $t = 0$, but a *random variable*. Consequently, the computation of $\mathrm{E}[\widetilde{B}_h^{(t)} | X_t]$ is not trivial since $\widetilde{B}_h^{(t)}$ and $X_t$ are not independent.

### 4.1. Proof of Theorem 2, part (2)

The proof of this part is actually a direct consequence of Theorem 2.4 in [7]. But, for the sake of completeness on one hand and taking into account that we are dealing with the one-dimensional case on the other hand, we give all the details here. Indeed, contrary to the multidimensional case, it is easy to compute explicitly the coefficients which appear (see Lemma 1 below and compare with [1] Theorem 31 or [7] Proposition 5.4) which also has its own interest from our point of view.

The differential operators $\Gamma_I$ appearing in Theorem 2 are recursively[*] defined by

$$\Gamma_{(0)}(f,b) = bf', \qquad \Gamma_{(1)}(f,b) = f'$$

and, for $I \in \{0,1\}^k$,

$$\Gamma_{(I,0)}(f,b) = b(\Gamma_I(f,b))', \qquad \Gamma_{(I,1)}(f,b) = (\Gamma_I(f,b))' \tag{17}$$

with $(I,0), (I,1) \in \{0,1\}^{k+1}$. The constants $c_I$ are explained as follows. Set

$$dB_t^{(i)} = \begin{cases} dB_t, & \text{if } i = 1, \\ dt, & \text{if } i = 0. \end{cases} \tag{18}$$

Then for a sequence

$$I = (i_1, \ldots, i_k) \in \{0,1\}^k,$$

we define

$$c_I = \mathrm{E}\left[\int_{\Delta^k[0,1]} dB^I\right] = \mathrm{E}\left[\int_0^1 dB_{t_k}^{(i_k)} \int_0^{t_k} dB_{t_{k-1}}^{(i_{k-1})} \cdots \int_0^{t_2} dB_{t_1}^{(i_1)}\right]. \tag{19}$$

Set $|I| = \sum_{1 \leq j \leq k} i_j$. Equivalently, $|I|$ denotes the number of integrals with respect to '$dB$'. Since $B$ and $-B$ have the same law, note that we have

$$c_I = c_I(-1)^{|I|}.$$

Thus, $c_I$ vanishes when $|I|$ is odd. In general, the computation of the coefficients $c_I$ can be made as follows.

---

[*]We can also use a rooted trees approach in order to define the $\Gamma_I$'s. See [7] for a thorough study, even in the multidimensional case and where $H > 1/3$.

**Lemma 1.** *Let $I \in \{0,1\}^k$. We denote by $J = \{j_1 < \cdots < j_m\}$ the set of indices $j \in \{1,\ldots,k\}$ such that $dB_t^{(j)} = dt$. We then have that $c_I$ is given by*

$$\int_0^1 dt_{j_m} \cdots \int_0^{t_{j_2}} dt_{j_1} \mathrm{E}\left[\frac{(B_1 - B_{t_{j_m}})^{k-j_m}(B_{t_{j_1}})^{j_1-1}}{(k-j_m)!(j_1-1)!} \prod_{k=2}^m \frac{(B_{t_{j_k}} - B_{t_{j_{k-1}}})^{j_k-j_{k-1}}}{(j_k - j_{k-1})!}\right].$$

The expectation appearing in the above formula can always be computed using the moment generating function of an $m$-dimensional Gaussian random variable. For instance, we have

$$\mathrm{E}\left[\int_0^1 \int_0^{t_3} \int_0^{t_2} dt_1\, dB_{t_2}\, dB_{t_3}\right] = \frac{1}{2(2H+1)},$$

$$\mathrm{E}\left[\int_0^1 \int_0^{t_3} \int_0^{t_2} dB_{t_1}\, dt_2\, dB_{t_3}\right] = \frac{2H-1}{2(2H+1)},$$

$$\mathrm{E}\left[\int_0^1 \int_0^{t_3} \int_0^{t_2} dB_{t_1}\, dB_{t_2}\, dt_3\right] = \frac{1}{2(2H+1)},$$

$$\mathrm{E}\left[\int_0^1 \int_0^{t_4} \int_0^{t_3} \int_0^{t_2} dt_1\, dt_2\, dB_{t_3}\, dB_{t_4}\right] = \frac{1}{2(2H+1)(2H+2)},$$

$$\mathrm{E}\left[\int_0^1 \int_0^{t_4} \int_0^{t_3} \int_0^{t_2} dt_1\, dB_{t_2}\, dt_3\, dB_{t_4}\right] = \frac{H}{(2H+1)(2H+2)}, \qquad (20)$$

$$\mathrm{E}\left[\int_0^1 \int_0^{t_4} \int_0^{t_3} \int_0^{t_2} dt_1\, dB_{t_2}\, dB_{t_3}\, dt_4\right] = \frac{1}{2(2H+1)(2H+2)},$$

$$\mathrm{E}\left[\int_0^1 \int_0^{t_4} \int_0^{t_3} \int_0^{t_2} dB_{t_1}\, dt_2\, dB_{t_3}\, dt_4\right] = \frac{H}{(2H+1)(2H+2)},$$

$$\mathrm{E}\left[\int_0^1 \int_0^{t_4} \int_0^{t_3} \int_0^{t_2} dB_{t_1}\, dt_2\, dt_3\, dB_{t_4}\right] = \frac{H(2H-1)}{2(2H+1)(2H+2)},$$

$$\mathrm{E}\left[\int_0^1 \int_0^{t_4} \int_0^{t_3} \int_0^{t_2} dB_{t_1}\, dB_{t_2}\, dt_3\, dt_4\right] = \frac{1}{2(2H+1)(2H+2)}.$$

**Lemma 2.** *When $f : \mathbb{R} \to \mathbb{R} \in C_b^\infty$, we have*

$$f(X_h) = f(x) + \sum_{k=1}^{n-1} \sum_{I_k \in \{0,1\}^k} \Gamma_{I_k}(f,b)(x) \int_{\Delta^k[0,h]} dB^{I_k}(t_1,\ldots,t_k)$$

$$+ \sum_{I_n \in \{0,1\}^n} \int_{\Delta^n[0,h]} \Gamma_{I_n}(f,b)(X_{t_1}) dB^{I_n}(t_1,\ldots,t_n), \qquad (21)$$

where, again using the convention (18), for $g:\mathbb{R} \to \mathbb{R} \in C_b^\infty$,

$$\int_{\Delta^k[0,h]} g(X_{t_1})\,dB^{I_k}(t_1,\ldots,t_k) \triangleq \int_0^h dB_{t_k}^{(i_k)} \int_0^{t_k} dB_{t_{k-1}}^{(i_{k-1})} \cdots \int_0^{t_2} dB_{t_1}^{(i_1)} g(X_{t_1}).$$

**Proof.** Applying (3) twice, we can write

$$f(X_h) = f(x) + \int_0^h f'(X_s)\,dB_s + \int_0^h (bf')(X_s)\,ds$$
$$= f(x) + \Gamma_{(1)}(f,b)(x)\,B_h + \Gamma_{(0)}(f,b)(x)h$$
$$+ \sum_{I_2 \in \{0,1\}^2} \int_{\Delta^2[0,h]} \Gamma_{I_2}(f,b)(X_{t_1})\,dB^{I_2}(t_1,t_2).$$

Applying (3) repeatedly, we finally obtain (21). $\square$

The remainder can be bounded by the following lemma.

**Lemma 3.** *If $n \geq 2$, $\varepsilon > 0$ (small enough) and $g:\mathbb{R} \to \mathbb{R} \in C_b^\infty$ are fixed, we have*

$$\sum_{I_n \in \{0,1\}^n} E\left|\int_{\Delta^n[0,h]} g(X_{t_1})\,dB^{I_n}(t_1,\ldots,t_n)\right| = O(h^{nH-\varepsilon}).$$

**Proof.** This involves a direct application of Theorem 2.2 in [6], combined with the Garsia, Rodemich and Rumsey Lemma [4]. $\square$

Thus, in order to obtain the asymptotic expansion of $\widetilde{P}_t f(h)$, Lemmas 2 and 3 say that it is sufficient to compute

$$E\left[\int_{\Delta^k[0,h]} dB^{I_k}(t_1,\ldots,t_k)\right]$$

for any $I_k \in \{0,1\}^k$, with $1 \leq k \leq n-1$. By the self-similarity and the stationarity of fractional Brownian motion, we have that

$$\int_{\Delta^k[0,h]} dB^{I_k}(t_1,\ldots,t_k) \stackrel{\mathcal{L}}{=} h^{H|I_k|+k-|I_k|} \int_{\Delta^k[0,1]} dB^{I_k}(t_1,\ldots,t_k).$$

Hence, it follows that

$$E\left[\int_{\Delta^k[0,h]} dB^{I_k}(t_1,\ldots,t_k)\right] = h^{H|I_k|+k-|I_k|} c_{I_k}$$

and the proof of point (2) of Theorem 2 is a consequence of Lemmas 2 and 3 above. $\square$

### 4.2. Proof of Theorem 2, point (1)

Let $f:\mathbb{R} \to \mathbb{R} \in C_b^\infty$ and $X$ be the solution of (5). We then know (see Section 2) that $X$ is a fractional Brownian motion with Hurst index $H$ under the new probability $\mathbf{Q}$ defined by $d\mathbf{Q} = \eta^{-1} d\mathbf{P}$, with $\eta$ given by (10). Since $b:\mathbb{R} \to \mathbb{R} \in C_b^\infty$, we observe that $\eta \in \mathbb{D}^{1,2}$. Moreover, the following well-known formula holds for any $\xi \in L^2(\mathbf{P}) \cap L^2(\mathbf{Q})$:

$$\mathrm{E}[\xi|X_t] = \frac{\mathrm{E}^{\mathbf{Q}}[\eta \xi|X_t]}{\mathrm{E}^{\mathbf{Q}}[\eta|X_t]}.$$

In particular,

$$\mathrm{E}[f(X_{t+h}) - f(X_t)|X_t] = \frac{\mathrm{E}^{\mathbf{Q}}[\eta(f(X_{t+h}) - f(X_t))|X_t]}{\mathrm{E}^{\mathbf{Q}}[\eta|X_t]}.$$

We now need the following technical lemma.

**Lemma 4.** *Let $\zeta \in \mathbb{D}^{1,2}(\mathfrak{H})$ be a random variable. Then for any $h > 0$, the conditional expectation $\mathrm{E}[\zeta(f(B_{t+h}) - f(B_t))|B_t]$ is equal to*

$$H(2H-1)f'(B_t) \int_0^T du\, \mathrm{E}[D_u\zeta|B_t] \int_t^{t+h} |v-u|^{2H-2} dv$$

$$+ \frac{1}{2} t^{-2H} f'(B_t)(B_t \mathrm{E}[\zeta|B_t] - \mathrm{E}[\langle D\zeta, \mathbf{1}_{[0,t]}\rangle_{\mathfrak{H}}|B_t])(h^{2H} - (t+h)^{2H} + t^{2H})$$

$$- \frac{H}{2} t^{-2H} f''(B_t) \mathrm{E}[\zeta|B_t]$$

$$\times \int_t^{t+h} ((v-t)^{2H-1} - v^{2H-1})(t^{2H} + v^{2H} - (v-t)^{2H}) dv$$

$$+ \frac{1}{2} f''(B_t) \mathrm{E}[\zeta|B_t]((t+h)^{2H} - t^{2H})$$

$$+ H(2H-1) \int_0^T du \int_t^{t+h} |v-u|^{2H-2} \mathrm{E}[D_u\zeta(f'(B_v) - f'(B_t))|B_t] dv$$

$$+ Ht^{-2H} B_t \int_t^{t+h} ((v-t)^{2H-1} - v^{2H-1}) \mathrm{E}[\zeta(f'(B_v) - f'(B_t))|B_t] dv \qquad (22)$$

$$- Ht^{-2H} \int_t^{t+h} ((v-t)^{2H-1} - v^{2H-1}) \mathrm{E}[\langle D\zeta, \mathbf{1}_{[0,t]}\rangle_{\mathfrak{H}}(f'(B_v) - f'(B_t))|B_t] dv$$

$$- \frac{H}{2} t^{-2H} \int_t^{t+h} ((v-t)^{2H-1} - v^{2H-1})(t^{2H} + v^{2H} - (v-t)^{2H})$$

$$\times \mathrm{E}[\zeta(f''(B_v) - f''(B_t))|B_t] dv$$

$$+ H \int_t^{t+h} \mathrm{E}[\zeta(f''(B_v) - f''(B_t))|B_t]v^{2H-1}\,\mathrm{d}v.$$

**Proof.** Let $g: \mathbb{R} \to \mathbb{R} \in C_b^1$. We can write, using the Itô formula ([9], (5.44), page 294) and basics identities of Malliavin calculus,

$$\mathrm{E}[\zeta g(B_t)(f(B_{t+h}) - f(B_t))]$$
$$= \mathrm{E}[\zeta g(B_t)\delta(f'(B.)\mathbf{1}_{[t,t+h]})] + H\int_t^{t+h} \mathrm{E}[\zeta g(B_t)f''(B_v)]v^{2H-1}\,\mathrm{d}v$$
$$= \mathrm{E}[g(B_t)\langle D\zeta, \mathbf{1}_{[t,t+h]}f'(B.)\rangle_{\mathfrak{H}}] + \mathrm{E}[\zeta g'(B_t)\langle \mathbf{1}_{[0,t]}, \mathbf{1}_{[t,t+h]}f'(B.)\rangle_{\mathfrak{H}}]$$
$$\quad + H\int_t^{t+h} \mathrm{E}[\zeta g(B_t)f''(B_v)]v^{2H-1}\,\mathrm{d}v$$
$$= H(2H-1)\int_0^T \mathrm{d}u\int_t^{t+h} |v-u|^{2H-2}\mathrm{E}[f'(B_v)g(B_t)D_u\zeta]\,\mathrm{d}v$$
$$\quad + H(2H-1)\int_0^t \mathrm{d}u\int_t^{t+h} (v-u)^{2H-2}\mathrm{E}[\zeta f'(B_v)g'(B_t)]\,\mathrm{d}v$$
$$\quad + H\int_t^{t+h} \mathrm{E}[\zeta g(B_t)f''(B_v)]v^{2H-1}\,\mathrm{d}v.$$

But,

$$\mathrm{E}[\zeta g(B_t)f'(B_v)B_t] = \mathrm{E}[\langle D(\zeta g(B_t)f'(B_v)), \mathbf{1}_{[0,t]}\rangle_{\mathfrak{H}}]$$
$$= \mathrm{E}[g(B_t)f'(B_v)\langle D\zeta, \mathbf{1}_{[0,t]}\rangle_{\mathfrak{H}}] + \mathrm{E}[\zeta g'(B_t)f'(B_v)]t^{2H}$$
$$\quad + \tfrac{1}{2}\mathrm{E}[\zeta g(B_t)f''(B_v)](t^{2H} + v^{2H} - (v-t)^{2H}).$$

Consequently,

$$\mathrm{E}[\zeta g(B_t)(f(B_{t+h}) - f(B_t))]$$
$$= H(2H-1)\int_0^T \mathrm{d}u\int_t^{t+h} |v-u|^{2H-2}\mathrm{E}[f'(B_v)g(B_t)D_u\zeta]\,\mathrm{d}v$$
$$\quad + H(2H-1)t^{-2H}\int_0^t \mathrm{d}u\int_t^{t+h} (v-u)^{2H-2}\mathrm{E}[\zeta g(B_t)f'(B_v)B_t]\,\mathrm{d}v$$
$$\quad - H(2H-1)t^{-2H}\int_0^t \mathrm{d}u\int_t^{t+h} (v-u)^{2H-2}\mathrm{E}[g(B_t)f'(B_v)\langle D\zeta, \mathbf{1}_{[0,t]}\rangle_{\mathfrak{H}}]\,\mathrm{d}v$$
$$\quad - \frac{1}{2}H(2H-1)t^{-2H}\int_0^t \mathrm{d}u\int_t^{t+h} (v-u)^{2H-2}\mathrm{E}[\zeta g(B_t)f''(B_v)]$$
$$\qquad\qquad \times (t^{2H} + v^{2H} - (v-t)^{2H})\,\mathrm{d}v$$

$$+ H \int_t^{t+h} \mathrm{E}[\zeta g(B_t) f''(B_v)] v^{2H-1} \, \mathrm{d}v.$$

We deduce that

$$\mathrm{E}[\zeta (f(B_{t+h}) - f(B_t)) | B_t]$$
$$= H(2H-1) \int_0^T \mathrm{d}u \int_t^{t+h} |v-u|^{2H-2} \mathrm{E}[f'(B_v) D_u \zeta | B_t] \, \mathrm{d}v$$
$$+ H t^{-2H} B_t \int_t^{t+h} ((v-t)^{2H-1} - v^{2H-1}) \mathrm{E}[\zeta f'(B_v) | B_t] \, \mathrm{d}v$$
$$- H t^{-2H} \int_t^{t+h} ((v-t)^{2H-1} - v^{2H-1}) \mathrm{E}[f'(B_v) \langle D\zeta, \mathbf{1}_{[0,t]} \rangle_{\mathfrak{H}} | B_t] \, \mathrm{d}v$$
$$- \frac{H}{2} t^{-2H} \int_t^{t+h} ((v-t)^{2H-1} - v^{2H-1})(t^{2H} + v^{2H} - (v-t)^{2H}) \mathrm{E}[\zeta f''(B_v) | B_t] \, \mathrm{d}v$$
$$+ H \int_t^{t+h} \mathrm{E}[\zeta f''(B_v) | B_t] v^{2H-1} \, \mathrm{d}v.$$

Finally, (22) follows. □

First, we apply the previous lemma with $\zeta = \eta$, $\mathrm{E} = \mathrm{E}^{\mathbf{Q}}$ and $B = X$, with $\eta$ given by (10), $\mathrm{d}\mathbf{Q} = \eta^{-1} \mathrm{d}\mathbf{P}$ and $X$ given by (5). We note that $\eta \in \mathbb{D}^{\infty,2}$ (see, e.g., [9] Lemma 6.3.1 and [2] for the expression of Malliavin derivatives via the transfer principle). We can particularly deduce that each random variable $V_k$, recursively defined by $V_0 = \eta$ and $V_{k+1} = \langle DV_k, \mathbf{1}_{[0,t]} \rangle_{\mathfrak{H}}$ for $k \geq 0$, belongs to $\mathbb{D}^{1,2}$.

In (22), the deterministic terms $\int_t^{t+h} |v-u|^{2H-2} \, \mathrm{d}v$, $(t+h)^{2H} - t^{2H}$ and

$$\int_t^{t+h} ((v-t)^{2H-1} - v^{2H-1})(t^{2H} + v^{2H} - (v-t)^{2H}) \, \mathrm{d}v$$

have a Taylor expansion in $h$ of the type (15). Lemma 4 allows the first term of the asymptotic expansion to be obtained using the fact that $\int_t^{t+h} \phi(s) \, \mathrm{d}s = h\phi(t) + o(h)$ for any continuous function $\phi$. By a recursive argument, again using Lemma 4, we finally deduce that (15) holds. This concludes the proof of Theorem 2. □

## Acknowledgements

The computations (20) were carried at by A. Neuenkirch. We would like to thank him.